\documentclass{osudgm}
\usepackage{amssymb,amsmath,cmmib57}

\def\cite#1{{\bf[#1]}} 

\newtheorem{theorem}{Theorem}[section] 

\newtheorem{corollary}[theorem]{Corollary}

\theoremtextfont{\rmfamily}

\title{Balanced presentations of the trivial group on two generators and  the  Andrews-Curtis conjecture}

\author{Alexei D. Miasnikov and Alexei G. Myasnikov}
\begin{document}
\maketitle

\begin{abstract}
The Andrews-Curtis conjecture states  that every balanced
presentation of the trivial group can be reduced to the standard one by a
sequence of the elementary Nielsen transformations and conjugations. In this
paper we describe all balanced presentations of the trivial group on two
generators and with the total length of relators $\leq 12$. We show that all
these presentations satisfy the Andrews-Curtis conjecture.
\end{abstract}

\begin{classification}
Primary 20E05, 20F05, 68T05; Secondary  57M05,
57M20. \\  This research was partially supported by the NSF, grant DMS-9973233.
\end{classification}

\section{ Introduction }

  Let  $F = F(X)$   be a free group of  rank
$n \geq 2$ with a basis $X = {\{}x_1, ..., x_n {\}}$. Consider the following
transformations of an $n$-tuple $W = (w_1, \ldots, w_n)$
 of elements from  $F$:
\smallskip

\noindent {\bf (AC1)} replace $w_i$   by  $w_i w_j,$  \  $j \ne i$;
\smallskip

\noindent {\bf (AC2)} replace $w_i$  by $w_i^{-1}$;
\smallskip

\noindent {\bf (AC3)} replace $w_i$ by  $f w_i f^{-1}$ for some
$f \in F$,
\smallskip

\noindent and leave $w_k$ fixed for all $k \neq i.$
\smallskip

\noindent The transformations (AC1) and (AC2) are usually called  elementary Nielsen
transformations. We will refer to the transformations (AC1) - (AC3) as the
{\em elementary AC-transformations}.

  Let $F^n$  be the cartesian product of $n$ copies of the group $F.$ Two $n$-tuples
   $V$ and $W $  from $F^n$   are called
{\em Andrews-Curtis equivalent} (or {\em AC-equivalent}) if one of them can be
obtained from the other by a finite sequence of elementary AC-transformations.
 In this event,
 we write $V \sim W$. The relation $\sim $ is an equivalence  relation
 on the set $F^n.$ The following conjecture, which appears to
be of interest in topology as well as group theory,  was raised by
J.J. Andrews and M.L. Curtis in \cite{AC}.

\smallskip
\noindent {\bf The Andrews-Curtis conjecture.} Elements $w_1, \ldots , w_n \in F(X)$
generate $F(X)$ as a normal subgroup if and only if $(w_1, \ldots ,w_n)
\sim (x_1, \ldots , x_n)$.
\smallskip

One can formulate this conjecture in terms of presentations.
 A group presentation $\langle x_1, \ldots, x_m \mid r_1, \ldots,r_n
\rangle$ is called {\em balanced} if $m = n.$  The balanced presentation  $\langle x_1, \ldots, x_n \mid x_1, \ldots,x_n
\rangle$ of the trivial group is called {\em standard}.  We shall say that two
presentations with the same generators are AC-equivalent if the tuples of
relators in these presentations are AC-equivalent. Plainly, the AC-conjecture
is equivalent to the following one: {\em every balanced presentation of the
trivial group is AC-equivalent to the standard one with the same set of
generators}.

 There is a survey \cite{BM} by
R.G. Burns and O. Macedonska on the AC-conjecture from group theory viewpoint. For
relevant topological results we refer to a survey \cite{HM} by C. Hog-Angeloni
and W. Metzler. The prevailing opinion seems to be that the AC-conjecture is
false. Moreover, several potential counterexamples are known. We say that a
balanced presentation of the trivial group $\langle x_1,
\ldots, x_n \mid w_1,
\ldots , w_n \rangle$ is a {\em potential counterexample} to the AC-conjecture
if, firstly, it is not known to be AC-equivalent to the standard one, and,
secondly, no one of the elementary AC-transformations decreases the total
length $|w_1| + \ldots + |w_n|$ of the relators. Below we list some of the
shortest and the most established potential counterexamples:

\smallskip
\noindent {\bf (1)} $\langle x,y| ~x^{-1}y^2 x = y^3, ~y^{-1}x^2 y = x^3\rangle$.

\smallskip
\noindent {\bf (2)} $\langle x,y, z| ~y^{-1}x y = x^2, ~z^{-1}y z = y^2,
~x^{-1}z x = z^2\rangle$.

\smallskip
\noindent {\bf (3)} $\langle x,y| ~x^3 = y^4, ~xyx=yxy\rangle$.

\smallskip
\noindent The first two presentations have been known for almost 20 years; for a
discussion we refer to the survey \cite{BM}. Example (3) is the second
presentation in the series

\smallskip
\noindent {\bf (4)} $\langle x,y| ~x^n = y^{n+1}, ~xyx=yxy\rangle,  ~n \ge 2$,

\smallskip
\noindent which is due to S. Akbulut and R. Kirby \cite{AK}. This series  has been known for
15 years. Notice that for $n = 2$ the series (4) gives the presentation

\smallskip
\noindent {\bf (5)} $\langle x,y| ~x^2 = y^3, ~xyx=yxy\rangle, $

\smallskip
\noindent  which, until recently, has also been considered as a potential
counterexample. In \cite{MA} the first  author proved that this
presentation satisfies the AC-conjecture. For this purpose he designed a
genetic algorithm to search for corresponding sequences of the elementary
AC-transformations (see \cite{MA} for details). This algorithm is also
instrumental for the following theorem which is the main result of this paper.

\begin{theorem}
\label{th:1}
Every balanced presentation of the trivial group on two generators with the
total length of the relators at most 12 satisfies the Andrews-Curtis conjecture.
\end{theorem}

The idea of the proof of this theorem is very simple: we list all  balanced
presentations of the trivial group on two generators where  the total length of
relators is at most 12:
$$
{\cal P} = \{P_1, P_2, P_3, \ldots \}
$$
and then apply the genetic algorithm to each of the listed presentations.
However, there are two issues to address here. The first one concerns  the
listing of the trivial presentations. Indeed, it is known that there is no
algorithm to decide whether a group given by a finite presentation is trivial
or not \cite{Adjan}, \cite{Rabin}. Whether  such an algorithm exists for balanced
presentations is an open and difficult problem \cite{BMS}. In the particular
case when the total length of relators in the presentations is at most 12, the
following result enables one to list the set ${\cal P}.$

\begin{theorem} Let $G$ be a group defined by a presentation
\label{th:2}
\begin{description}
\item {\bf (6)}  $\langle \ x,y \mid r(x,y),\  s(x,y) \ \rangle$, where
$|r(x,y)| + |s(x,y)| \leq 12$.
\end{description}
 If the abelianization of $G$ is trivial then $G$ is either the trivial
group or $G$ is isomorphic to the following finite group of order 120
$$
\langle x,y \mid y x y = x^2, \ x y x = y^4\  \rangle.
$$
\end{theorem}

\noindent    The proof of this theorem is based on computer computations with the
software package Magnus. We discuss the proof in Section 2.

The second issue is related to the real time required to carry out the
computations with the genetic algorithm. It turns out that there are  about $10^6$
presentations in the set ${\cal P}.$ To run the genetic algorithm on each of
them would take too much time. So one has to exploit tricks and shortcuts
(pribambases) to decrease the time. We discuss this in Section 2. However, it
is worthwhile to note here that the presentations from ${\cal P}$ with
total length up to 10 are relatively easy to reduce to the standard one by
AC-transformations. So, the minimal total length of relators in non-trivial
examples from ${\cal P}$ appears to be 11. Surprisingly enough, all the
"difficult" presentations from ${\cal P}$ with  total length 11 (which
cannot be easily reduced to the standard one) are readily seen to be
AC-equivalent to the presentation {\bf (5)}: $\langle x,y| ~x^2 = y^3,
~xyx=yxy\rangle.$ The most difficult examples to ``crack" are those with
total length 12, which cannot be easily reduced neither to the standard one
nor to the presentation {\bf (5)} (of length 11). We list these most
interesting presentations in the corollary below.

\begin{corollary} The  following presentations
\label{cor:1}
of the trivial group are AC-equivalent to $\langle x,y| ~x, y\rangle$:

$$
\langle x,y| ~x^{-1}y^2 x = y^3, ~x^2=y^\epsilon x y^\delta\rangle,
$$
  where $\epsilon, \delta \in \{1,-1\}.$
\end{corollary}

  One can find in \cite{MA} the corresponding sequences of
AC-transformations that reduce the presentations above to the standard
presentation of the trivial group.

  It follows now from Theorem \ref{th:1} that  the minimal total length of relators in
  potential counterexamples, that still stand, is 13. This allows us to formulate the
  following

\begin{corollary}
\label{cor:2}
 The presentation {\bf (3)}:
$$\langle x,y | ~x^3 = y^4, ~xyx=yxy\rangle$$
is, at the present time, a minimal potential counterexample to the
Andrews-Curtis conjecture.
\end{corollary}

\section{Description of the algorithms }
\label{sec:results}

We start with few known algorithms, which play an important part in our
proofs.

In 1936 J.H.C. Whitehead \cite{Whit} gave an algorithm which for given $m$-tuples
$U=(u_1,...,u_m)$ and $V=(v_1,...,v_m)$ from $F^m$ decides whether there
exists an automorphism $\phi \in Aut(F)$ such that $\phi(u_i) = v_i, i
= 1, \ldots,m.$ Furthermore, if  such an automorphism  exists, then
the algorithm finds one. In general, the time-complexity of the Whitehead
algorithm is exponential. However, in the particular case when one needs to
check whether a given element $f \in F$ can be mapped by an automorphism of
$F$ to the element $x_1,$ the Whitehead method is polynomial in time with
respect to the length of $f$ (for a given fixed $F $). Recall that an element
$f \in F$ is called {\em primitive} if $\phi(f) = x_1$ for some $\phi \in
Aut(F).$ It follows that we can recognize primitive elements in a free group
of rank two quite effectively.

The other algorithm that we used in our proofs is the Todd-Coxeter algorithm
\cite{TC}, \cite{Jon}. This is a systematic procedure for enumerating cosets
of a given finitely generated subgroup of a given finitely presented group. In
particular, if the group given by a finite presentation is finite, then  the
Todd-Coxeter algorithm will eventually recognize this, and it will give a
multiplication table of the group. At the present time, the Todd-Coxeter
algorithm (and its variations) are among the most powerful application of
computers to  group theory (see \cite{Havas} for details).

The third algorithm we want to mention here allows one to check whether the
abelianization $G/[G,G]$ of the group $G$ given by a finite presentation
$\langle x_1, \ldots , x_n \mid r_1, \ldots , r_m \rangle$ is trivial or not.

This algorithm is relatively fast (at least in our case); it calculates the
canonical invariants of the abelian group $G/[G,G]$. One can find a complete
description of the algorithm in the book \cite{Sims} by C.Sims.

Now we explain how these methods can be used in proving Theorems \ref{th:1} and \ref{th:2}. We
combine both proofs into a single procedure that was carried out by a
computer. This procedure consists of the following steps.

\begin{enumerate}
\item Generate a list $L_1$ of balanced presentations on  two generators
and with the total length of relators $\leq 12$. There are about $6\cdot 10^6$
of such presentations.

\item For each presentation $P$ in the list $L_1$,  check to see whether the
abelianization of the group defined by $P$ is trivial or not. Delete, one by
one, all presentations from $L_1$ with a non-trivial abelianization. Denote by
$L_2$ the resulting list. Plainly, all the trivial groups from $L_1$ are also
in $L_2.$ There are about $10^6$ presentations in $L_2.$

\item  Apply the Whitehead algorithm to each of the relators in
every presentation $P$ in $L_2$ to check whether the relator is a primitive
element in $F$ or not. If the relator is primitive then the group defined by
$P$ is trivial. Moreover, it is not hard to see that, in this event, the
presentation $P$ is AC-equivalent to the standard one \cite{MA}. Delete, one
by one, all the presentations from $L_2$ with primitive elements among its
relators. All the deleted presentations satisfy the AC-conjecture. Denote by
$L_3$ the resulting list. There are 122240 presentations in $L_3$.

\item  Observe that a cyclic permutation of a relator is a particular case of
the transformation $(AC3)$. We can use this to reduce the number of
presentations in $L_3$ which belong to the same AC-equivalence class. Compare
presentations in $L_3$, by cyclically permuting their relations, and leave
only one presentation from each equivalence class. Denote the resulting list
by $L_4.$ Only 1648 presentations are left in the list $L_4$.

\item  At this point we want to sort out the presentations in $L_4$ which define
the trivial group. We apply the Todd-Coxeter algorithm to each presentation
$P$ from $L_4$ to compute the order of the group defined by $P.$ Notice that
the Todd-Coxeter algorithm is time-consuming, so we want to apply it to as few
presentations as possible. Luckily, all the groups in the list $L_4$ happened
to be finite, so the Todd-Coxeter algorithm eventually stopped and gave the
answer. It turns out that all the groups in $L_4$ are trivial, except 16
groups of order 120. It follows now that all the groups from $L_2$ are either
trivial or of order 120. To finish the proof of Theorem \ref{th:2} it suffices to
notice that all these groups are isomorphic to each other, and hence they are
isomorphic to this particular one, given by the presentation
$$
\langle x,y \mid y x y = x^2, \ x y x = y^4\ \rangle.
$$
Indeed, every group of  order 120 with  trivial abelianization is a
central extension of a cyclic group of order 2 by the simple group $A_5$. All
such groups are isomorphic. Denote by $L_5$ the list of all presentations from
$L_4$ which define the trivial group. There are 1632 presentations in $L_5.$

\item  This is the last and the most time-consuming  step. We apply the genetic
algorithm from \cite{MA} (and some of its variations) to check whether the
presentations from $L_5$ satisfy the AC-conjecture. We show that every
presentation from $L_5$ is AC-equivalent either to the standard presentation
of the trivial group, or to the presentation {\bf (5)}, or to one of the
presentations in Corollary \ref{cor:1}, which are already known to satisfy the
AC-conjecture (see \cite{MA}).

At this step, all the presentations of the trivial group on two generators and
with  total length at most 12 have been shown to satisfy the AC-conjecture.
This proves Theorem \ref{th:1}.
\end{enumerate}

All  routines and procedures which we used here are available via the
Internet at \ {\em www.grouptheory.org}\ as a part of the software package
Magnus.

 Notice that in order to show that two presentations are AC-equivalent we
 use a modification of the genetic algorithm from  \cite{MA} in which  the fitness
 function is replaced by a new one. Namely, in this case we used the sum of
 the Hamming distances between  the relators as cyclic words. In  most of the
 occasions
 the genetic algorithm with this  new fitness  function worked  very fast. We refer
 to \cite{Hol} and \cite{Mit} for a general discussion on genetic algorithms.

\Coordinates
[Department of Mathematics, The City College of The City University of New York]\hfil\break
[paper mail address: Math. Dept., Convent av. at 138 st., The City College of CUNY, New York, NY, 10031]\hfil\break
[e-mail: alexei@groups.sci.ccny.cuny.edu]\hfil\break
\hfil\break
\endCoordinates


\frenchspacing
\begin{thebibliography}{A}
\bibitem[Adjan]{A} S.I. Adjan,On algorithmic problems in effectively 
complete classes of groups, Doklady Akad. Nauk SSSR  {\bf123} (1958), 13--16.

\bibitem[AK]{AK} S. Akbulut and R. Kirby, 
A potential smooth counterexample in dimension 4 to the Poincare conjecture, 
the Schoenflies conjecture, and the Andrews-Curtis conjecture, 
Topology \textbf{24} (1985), 375--390.
\bibitem[AC]{AC}J.J. Andrews and M.L. Curtis, Free groups and handlebodies,
Proc. Amer. Math. Soc. \textbf{16} (1965), 192--195.
\bibitem[BMS]{BMS} G. Baumslag, A. Myasnikov, and V. Shpilrain, 
Open problems in combinatorial group theory, in Groups, Languages and Geometry
(ed. by R.Gilman), Contemporary Math. {\bf250} (1999), 1--27.

\bibitem[BM]{BM}R.G. Burns and Olga Macedonska, Balanced presentations of
the trivial group, Bulletin London Math. Soc. \textbf{25} (1993), 513--526.

\bibitem[Havas]{Havas} John J. Cannon, Lucien A. Dimino, George Havas and
Jane M. Watson, Implementation and analysys of the
Todd-Coxeter algorithm, Math. Comput., \textbf{27}, 463--490.

\bibitem[HM]{HM} C. Hog-Angeloni, W. Metzler, The Andrews-Curtis conjecture and its
generalizations, in Two-dimensional homotopy and combinatorial group theory,
 London Math. Soc. Lecture Note Ser. {\bf197}, Cambridge Univ. Press,
Cambridge (1993), 365--380.

\bibitem[Hol]{Hol} J.H. Holland, Adaptation in Natural and Artificial
Systems, University of Michigan Press  (1975).
\bibitem[Jon]{Jon} D.L. Johnson, Topics in the theory of group
presentations, London Math. Soc. Lecture Note Ser. \textbf{42},
Cambridge-University press (1980).

\bibitem[MA]{MA} A.D. Miasnikov, Genetic algorithms and the  Andrews-Curtis
conjecture, The International Journal of Algebra and Computation vol.
\textbf{9} no. \textbf{6} (1999).

\bibitem[Mit]{Mit}M. Mitchell, An introduction to genetic algorithms,
Cambridge, MA:MIT Press (1996).

\bibitem[Ne]{Ne}B.H. Neumann,On some finite groups with trivial
multiplicator, Publ. Math. Debrecen \textbf{4} (1956), 190--194.

\bibitem[Rabin]{Rabin} M. O. Rabin, Recursive unsolvability of group theoretic
problems, Ann of Math. {\bf67} (1958), 172--194.

\bibitem[TC]{TC} J. A. Todd, H. S. M. Coxeter, A practical
method for enumerating cosets of a finite abstract group, Proc.
Edinburg Math. Soc.  \textbf{5}, (1936),  26--34.

\bibitem[Sims]{Sims} C. Sims, Computation With Finitely Presented
Groups, Cambridge Univ Press (1994).

\bibitem[Whit]{Whit}J. H. C. Whitehead, On equivalent sets of elements in a free
group, Ann. of Math. vol. \textbf{37} no. \textbf{4} (1936).





\end{thebibliography}
\end{document}